\documentclass[12pt,reqno]{amsart}

\usepackage{hyperref}

\newtheorem{Theorem}{Theorem}

\newtheorem{Conjecture}[Theorem]{Conjecture}
\theoremstyle{remark}

\allowdisplaybreaks[1]

\author[Michael J.\ Schlosser]{Michael J.\ Schlosser$^*$}
\address{Fakult\"at f\"ur Mathematik, Universit\"at Wien,
Oskar-Morgenstern-Platz~1, A-1090 Vienna, Austria}
\email{michael.schlosser@univie.ac.at}
\thanks{$^*$The author's research is partly supported by
FWF Austrian Science Fund grant P 32305.}

\title[A tribute to Dick Askey]{A tribute to Dick Askey\\[10pt]
\scriptsize A small contribution to the (September 15, 2019)\\
Liber Amicorum Richard ``Dick'' Allen Askey}

\begin{document}


\maketitle


Dick Askey has always been very supportive to the mathematical
community in general, and in particular to his mathematical
family: the Special Functions community which he cared about and
steadily nurtured (like a Godfather, but completely gently and nonviolently!).
I wish to thank him for his support and advice over the years
(during which he several times gave me encouraging
feedback on my work, was writing recommendation letters when
being asked, was giving me specific valuable advice on old literature
to read---for instance, by Heinrich August Rothe~\cite{rothe}, and by
Ferdinand Schweins~\cite{schweins}---, and so on).

It is difficult to estimate the value or the impact of a single person in one's life.
Dick Askey has had a steady light impact on me (just like water which
is inconspicuously dripping but is ultimately creating a canyon by erosion).
He was already one of the great shots in Special Functions
when I started to study Mathematics at the University of Vienna,
so for me he has always been around (and I took that for granted,
I never knew anything else).
It was in June 1996 (I just had completed my PhD thesis, working under the
direction of Christian Krattenthaler on multivariate basic hypergeometric series)
when I first met Dick, namely at the ``Miniconference
on $q$-Series'' which Gaurav Bhatnagar and Stephen Milne organized at
the Ohio State University. It was at this meeting where I first personally experienced
Dick Askey in action promoting special functions and supporting young people.
There were several other occasions where I was lucky to attend the same conference
where Dick was, including a conference on Ramanujan in Mysore, India, in
December 2012, to single out a somewhat (for me) exotic place where we met as well.

I now want to turn to mathematics, first abstractly, then concrete.
On this occasion of compiling a contribution to the Richard Askey Liber Amicorum
I would like to offer an (at least small) piece of mathematics
to Dick as a gift, as a (small) token of appreciation.
Various keywords come to (my) mind while contemplating about
Dick Askey and his work. {\em Special functions} and {\em $q$-series}
definitely belong to the main keywords. The imperative
``{\em Read the masters!}\/'' certainly comes to one's mind as well.
Dick has always stressed the importance of reading the
work of the old masters (Euler, etc.) for a better understanding of mathematics
(or in the creative art of actually doing mathematics), and argued that this
would also have a strong impact on one's research ability, either
by discovering what the old masters already knew or just by learning from their
thought process.
Last but not least, {\em positivity problems} and {\em problems of analytic flavor}
are intimately connected to Dick Askey's work too.
A mathematical gift to Dick should ideally connect various of the keywords
just mentioned.  (I admit that I have left out ``{\em orthogonal polynomials}"
but allow me a certain degree of artistic freedom to justify my --albeit arbitrary--
thought process!)

Consider the well-known binomial theorem,
\begin{equation}\label{binthm}
(1+z)^n=\sum_{k\ge 0}\binom nk z^k,
\end{equation}
with the binomial coefficient defined by
\begin{equation*}
\binom nk=\frac{n(n-1)\cdots(n-k+1)}{k!}.
\end{equation*}
In \eqref{binthm}, the exponent $n$ is a priori a nonnegative integer.
Isaac Newton, one of the old masters, experimented with this identity, formally
replaced the exponent $n$ by some fraction such as $\frac 12$, etc., and showed
that the identity still holds when both sides make sense.
(Today we know that in the nonterminating case, i.e., when $n$ is not a
nonnegative integer, $z$ is required to satisfy the condition $|z|<1$,
unless $z$ is considered a formal power series variable.)
Newton was thus the first person to consider a ``fractional'' extension of the
binomial theorem. In fact, as we know today, the integer parameter $n$
in \eqref{binthm} can be replaced by any complex number $\alpha$; the identity
then still holds (provided $|z|<1$, or if $z$ is simply a formal variable).

I would like to dedicate an observation, at this moment actually still a conjecture,
to Dick, which can be regarded to be a {\em fractional} extension of limiting
cases of the First and Second Borwein Conjectures (cf.\ \cite{andrews}).

Let $q$ be a complex variable with $0<|q|<1$.
As usual, the $q$-shifted factorial is defined as $(a;q)_0=1$ and
\begin{equation*}
(a;q)_n=(1-a)(1-aq)\cdots(1-aq^{n-1}).
\end{equation*}
This also makes sense for $n=\infty$; then the product has infinitely many factors.
For convenience, we shall also define
\begin{equation*}
(a_1,\dots,a_m;q)_n=(a_1;q)_n\cdots(a_m;q)_n
\end{equation*}
(of which we will only use the $m=2$ case in this tribute).

The celebrated First Borwein Conjecture (made by Peter Borwein around 1990, see
\cite{andrews})
asserts that for each nonnegative integer $n$,
the polynomials $A_n(q)$, $B_n(q)$ and $C_n(q)$ appearing in the dissection
\begin{equation}\label{borwein}
(q,q^2;q^3)_n=A_n(q^3)-qB_n(q^3)-q^2C_n(q^3)
\end{equation}
are actually polynomials in $q$ with {\em nonnegative} integer coefficients.

This conjecture was open for a long time and has only very recently been
settled by Chen Wang, in his 2019 doctoral thesis at the University of Vienna
(supervised by Christian Krattenthaler), see also \cite{wang} (which is part of his thesis).
Chen Wang's method of proof
(which follows a suggestion made by George Andrews in \cite{andrews})
is analytic in nature and makes careful use of asymptotic estimates to establish
bounds on the coefficients.

The First Borwein Conjecture is actually easy to show in the limit case $n=\infty$.
In that case one can make use of the famous Jacobi triple product identity to
prove the claimed nonnegativity of the series $A_\infty(q)$, $B_\infty(q)$
and $C_\infty(q)$.

The Second Borwein Conjecture (still being open) concerns a similar
dissection in terms of powers of $q$,
but with the expression on the left-hand side of \eqref{borwein} being squared.

Gaurav Bhatnagar (whom I know in person as long as Dick)
and I have recently formulated partial theta function
extensions of the first two Borwein Conjectures, see \cite{bhatschloss}.
There we replaced all the factors
in the respective $q$-shifted factorials by partial theta functions and observed
that similar positivity properties appear to hold.
My aim here is not to redeliver the results I have obtained with Gaurav
but rather to present something entirely new (and unspoiled!):

\begin{Conjecture}[Dedicated to Dick Askey]
Let $q$ be a complex number with $0<|q|<1$
(or view $q$ as a formal power series variable).
Further, let $d$ be a real number satisfying
\begin{equation*}
0.22799812734\ldots\approx\frac{9-\sqrt{73}}2\le d\le 1\quad\text{or}\quad
2\le d\le 3.
\end{equation*}
Then the series $A^{(d)}(q)$, $B^{(d)}(q)$, $C^{(d)}(q)$ appearing in
the dissection
\begin{equation}
(q,q^2;q^3)^d_\infty=A^{(d)}(q^3)-qB^{(d)}(q^3)-q^2C^{(d)}(q^3)
\end{equation}
are power series in $q$ with nonnegative real coefficients.
\end{Conjecture}

(The $d=3$ case of the above conjecture is the $n=\infty$ case of an 
observation made by Chen Wang, communicated to the community
by Christian Krattenthaler in his recent plenary talk
at OPSFA15 in Hagenberg, Austria, in July 2019.)  

\medskip
 Dick, I hope you like the above conjecture and can help to settle it!
 
 I wish you many more years of joy and productivity.
 Thanks a lot for serving as an idol and inspiration!

\end{document}